\documentclass{amsart}
\usepackage{amsmath,amsfonts,amssymb,amsthm,amsxtra,amscd,eucal}

\newcommand{\kk}{{\Bbbk}}

\newcommand{\zz}{\mathbb{Z}}

\DeclareMathOperator{\id}{id}

\newcommand{\cee}{{\mathcal{E}}}
\newcommand{\cff}{{\mathcal{F}}}
\newcommand{\ckk}{{\mathcal{K}}}

\newtheorem{thm}{Theorem}
\newtheorem{pr}[thm]{Proposition}
\newtheorem{corl}[thm]{Corollary}
\newtheorem{conj}[thm]{Conjecture}

\theoremstyle{definition}
\newtheorem{df}{Definition}
\newcommand{\thrm}[1]{\begin{thm}{#1}\end{thm}}
\newcommand{\prop}[1]{\begin{pr}{#1}\end{pr}}
\newcommand{\defn}[1]{\begin{df}{#1}\end{df}}
\newcommand{\cnl}[1]{\begin{corl}{#1}\end{corl}}
\newcommand{\cn}[1]{\begin{conj}{#1}\end{conj}}

\theoremstyle{remark}
\newtheorem*{rem}{Remark}

\newtheorem*{exm}{Example}
\newcommand{\exmp}[1]{\begin{exm}{#1}\end{exm}}
\newcommand{\prf}[1]{\begin{proof}{#1}\end{proof}}

\newcommand{\emp}{\varnothing}

\newcommand{\vphi}{\varphi}

\newcommand{\0}[1]{\overline{#1}}

\newcommand{\Ker}{\operatorname{Ker}}

\title[The category of $E_\infty$-coalgebras]
{The category of $E_\infty$-coalgebras, the $E_\infty$-coalgebra structure on the homology, and the dimension completion of the fundamental group}

\author{Grigory Rybnikov}

\address{National Research University Higher School of Economics}

\email{gr@mccme.ru, grybnikov@hse.ru}
\date{}
\begin{document}

\begin{abstract}
We study a special type of $E_\infty$-operads that govern strictly unital $E_\infty$-coalgebras (and algebras) over the ring of integers. Morphisms of coalgebras over such an operad are defined by using universal $E_\infty$-bimodules. Thus we obtain a category of $E_\infty$-coalgebras. It turns out that if the homology of an $E_\infty$-coalgebra have no torsion, then there is a natural way to define an $E_\infty$-coalgebra structure on the homology so that the resulting coalgebra be isomorphic to the initial $E_\infty$-coalgebra in our category. We also discuss some invariants of the $E_\infty$-coalgebra structure on homology and relate them to the invariant formerly used by the author to distinguish the fundamental groups of the complements of combinatorially equivalent complex hyperplane arrangements.
\end{abstract}

\maketitle

\section{Introduction}

Algebras and coalgebras over $E_\infty$-operads are usually discussed in relation to topology of simply connected topological spaces. For example, the relations between the homotopy category of pointed, simply connected CW-complexes and the category of $E_\infty$-coalgebras were studied by Justin Smith in \cite{Smith,Smith2}. Under certain assumptions, he showed that the cobar construction for an $E_\infty$-coalgebra can also be regarded as an $E_\infty$-coalgebra, and in the case of simply connected CW-complexes it is isomorphic (in the homotopy category of $E_\infty$-coalgebras) to the $E_\infty$-coalgebra on the chain complex of the loop space.

On the contrary, our main interest lies in the relations of $E_\infty$-coalgebras to the fundamental group. We construct a category of $E_\infty$-coalgebras so that there is a natural $E_\infty$-coalgebra structure on the homology of a topological space (provided that all homology groups are free Abelian), and the natural $E_\infty$-coalgebra structure on the chain complex of a topological space is isomorphic to the  $E_\infty$-coalgebra structure on its homology. Then, if we look at invariants of the $E_\infty$-coalgebra structure on homology with respect to isomorphisms in our category, we recognize among them (the dualized form of) Steenrod squares and Massey products. In particular, the invariant used in \cite{Ry} to distinguish the fundamental groups of the complements of combinatorially equivalent hyperplane arrangements (which is a form of dualized Massey product, see~\cite{Ry2}) turns out to be an invariant of the $E_\infty$-coalgebra structure on homology.

The author is deeply grateful to Anton Khoroshkin for very helpful discussion.

\section{Algebraic operads and modules over operads}

For information concerning algebraic operads, see the books \cite{KrizMay}, \cite{gr_loday_vallette}, \cite{MSS}, \cite{S2}, the article \cite{KM} and the preprint \cite{Smith2}. Following \cite{Smith2}, we will use the left action of symmetric groups rather than the right action, as in the majority of other sources. In calculations with complexes, we will follow the Koszul sign convention: when $a$ passes through $b$, we get the factor $(-1)^{\deg a\deg b}$. In particular, the tensor product of linear maps acts according to the formula $(f\otimes g)(x\otimes y)=(-1)^{\deg g\deg x}f(x)\otimes g(y)$.

Let $\kk$ be either a field or the ring of integers $\zz$, let $\ckk$ be the category of $\zz$-graded dg-$\kk$-modules with differential $\partial$ of degree $-1$ (i.e., the category of complexes of $\kk$-modules $\dots\leftarrow A_{-1}\leftarrow A_0\leftarrow A_1\leftarrow\dots$). The category $\ckk$ possesses the standard structure of symmetric monoidal category, with the direct sum as the coproduct and the tensor product as the monoidal product. The category of $\zz$-graded $\kk$-modules is embedded into $\ckk$ as the full subcategory consisting of complexes with zero differential; the category of $\kk$-modules can also be regarded as a full subcategory of $\ckk$ (consisting of complexes concentrated in degree $0$).

By $\Sigma$ denote the category with objects $[n]=\{1,\dots,n\}$, where $n$ runs through the set $\zz_+$ of nonnegative integers, and morphisms
$$
\Sigma([n],[m])=\begin{cases}\emp, & n\ne m,\\ \Sigma_n, & n=m,\end{cases}
$$
where $\Sigma_n$ is the symmetric group on the set $[n]$. Let $\ckk^\Sigma$ be the category of covariant functors $\Sigma\to\ckk$. The objects of this category will be called \emph{symmetric families of complexes}. For any $F\in\ckk^\Sigma$, we will write $F^{(n)}$ instead of $F([n])$, thus saying that we are given a symmetric family $F$ means saying that we are given a collection of complexes $F^{(n)}$ acted on by the corresponding symmetric groups. A morphism $f:A\to B$ of symmetric families is a collection of $\Sigma_n$-equivariant morphisms of complexes $\vphi^{(n)}:A^{(n)}\to B^{(n)}$.

Any symmetric family (a functor $F:\Sigma\to\ckk$) can be canonically extended to a functor to $\ckk$ from the category of all finite sets with bijections as morphisms. Denote this extension by the same letter $F$. We have
$$
F(S)=\Bigl(\bigoplus_{\vphi:S\to{[n]}}F^{(n)}\Bigr)_{\Sigma_n},
$$
where the direct sum is taken over the set of all bijections $\vphi:S\to{[n]}$, and $\Sigma_n$-coinvariants are taken with respect to the natural action of the symmetric group on this direct sum (the permutations are acting on the summands at the same time permuting them).

The category $\ckk^\Sigma$ possesses the following monoidal products: for $A,B\in\ckk^\Sigma$, we set
\begin{eqnarray}
(A\otimes B)^{(n)} &=& A^{(n)}\otimes B^{(n)}; \\
(A\boxtimes B)^{(n)} &=& \bigoplus_{[n]=I\sqcup J}A(I)\otimes B(J); \\
(A\odot B)^{(n)} &=& \bigoplus_{p\ge0}\left(\left(A^{\boxtimes p}\right)^{(n)}\mathbin{\mathop{\otimes}_{\Sigma_p}}B^{(p)}\right).
\end{eqnarray}
The first two of these products are symmetric, while the third one (which is often referred to as the \emph{plethysm}) is not. Each of these product provides $\ckk^\Sigma$ with a structure of monoidal category. We are interested in the third of them. The unit object for this structure is
$$
I^{(n)}_k=\begin{cases}\kk,& n=1,k=0,\\0,& n\ne1\text{\ or\ }k\ne0,\end{cases}
$$
where $k$ denotes the degree.

We can visualize the monoidal product $\odot$ (the plethysm) in the following way. Regard the elements of $A^{(n)}$ and $B^{(n)}$ as trees with a single vertex, one entry (from the right), and $n$ exits (to the left). The symmetric group $\Sigma_n$ acts on $A^{(n)}$ and on $B^{(n)}$ by permuting the exits. Then $A\odot B$ consists of trees of hight $2$ made as follows: the entry leads to a vertex of type $B$, and each exit of this vertex is connected to the entry of a vertex of type $A$. The symmetric group permutes all exits of this tree (taking account of the action of symmetric groups by permutations of the exits of the vertices of types $A$ and $B$).

Two morphisms $f,g:A\to B$ of the category $\ckk^\Sigma$ are said to be \emph{homotopic} if, for any $n\ge0$, there is a $\Sigma_n$-equivariant map $\vphi^{(n)}:A^{(n)}\to B^{(n)}$ (increasing the degree by $1$) such that $f^{(n)}-g^{(n)}=\vphi^{(n)}\circ\partial+\partial\circ\vphi^{(n)}$.

\defn{An \emph{operad} (or \emph{symmetric operad}) in the category of complexes is a monoid in the monoidal category $(\ckk^\Sigma,\odot)$. In other words, an operad is a symmetric family  $A$ equipped with an associative multiplication $\mu:A\odot A\to A$ and a unit element $1\in A(1)$ which satisfy natural compatibility conditions.

A left (right) \emph{module} over an operad $A$ is a symmetric family $B$ with a multiplications $\nu:A\odot B\to B$ ($\nu:B\odot A\to B$) such that the unit element acts as identity and natural associativity conditions are satisfied. A two-sided module will be referred to as a \emph{bimodule}.

Operads in the categories of $\kk$-modules and graded $\kk$-modules, as well as modules over them, are defined similarly.}

The operad structure can be described by means of the operation of ``$i$-th com\-po\-sition'' $\circ_i: A^{(m)}\otimes A^{(n)}\to A^{(m+n-1)}$, $i=1,\dots,n$, (substituting $A^{(m)}$ into $A(n)$ on $i$-th place). By definition, for $x\in A^{(m)}$ and $y\in A^{(n)}$,
$$
x\circ_i y=\mu^{(m+n-1)}((1\otimes\dots\otimes1\otimes x\otimes1\otimes\dots\otimes1)\otimes y),
$$
where $x$ is the $i$-th factor of $n$ in the tensor product $1\otimes\dots\otimes1\otimes x\otimes1\otimes\dots\otimes1$ (which belongs to the summand $A(\{1\})\otimes\dots\otimes A(\{i-1\})\otimes A(\{i,i+1,\dots,i+m-1\})\otimes A(\{i+m\})\otimes\dots\otimes A(\{n+m-1\})$ of the direct sum constituting $A^{\boxtimes n}$; we identify $A^{(1)}$ with $A(\{k\})$ and $A^{(m)}$ with $A(\{i,i+1,\dots,i+m-1\})$ by means of the monotone increasing bijections). A morphism of operads is a morphism of symmetric families which commutes with the composition operations $\circ_i$.

In calculations, for any collection $(a_1,\dots,a_n)$, $a_i\in A^{(m_i)}$, and any element $a\in A^{(n)}$, we will write $(a_1,\dots,a_n)\circ a$ rather than $\mu((a_1\otimes\dots\otimes a_n)\otimes a)$. Similarly, in the case of modules, the multiplication will also be denoted by the sign of composition~$\circ$.

For any complex $K$, the operads $\cee_K$ and $\cee^K$ are defined by the formulas $(\cee_K)^{(n)}=Hom(K^{{\otimes}n},K)$ and $(\cee^K)^{(n)}=Hom(K,K^{{\otimes}n})$ (here $Hom(K,L)$ denotes the complex with graded components $Hom_k(K,L)=\prod\limits_{l\in\zz}Hom(K_l,L_{l+k})$). In each case, the symmetric groups act by permutations of factors in the tensor product, the multiplication is given by composition, and the unit element is the identity map.

Similarly, for any two complexes $K,L$ we define the symmetric families $\cff_{K,L}$ and $\cff^{K,L}$ by the formulas $(\cff_{K,L})^{(n)}=Hom(K^{{\otimes}n},L)$ and $(\cff^{K,L})^{(n)}=Hom(K,L^{{\otimes}n})$. It is clear that $\cff_{K,L}$ has a natural structure of left module over $\cee_L$ and right module over $\cee_K$, while $\cff^{K,L}$ has a natural structure of left module over $\cee^K$ and right module over $\cee^L$.

Let $A$ be an operad, and let $K$ be a complex. We say that $K$ is given a structure of algebra (coalgebra) over $A$, if there is given a homomorphism of operads $A\to\cee_K$ (respectively, $A\to\cee^K$).

For example, consider the operad  ${Com}$, for which each component $Com^{(n)}$ is a rank $1$ free $\kk$-module $\kk e^{(n)}$ with generator $e^{(n)}$ of degree $0$ (the symmetric groups acts trivially); the composition is given by $e^{(m)}\circ_ie^{(n)}=e^{(m+n-1)}$ for all $m,n,i\le n$. Clearly, we have $e^{(1)}=1\in Com^{(1)}$. An algebra over this operad is just a commutative associative differential graded $\kk$-algebra with unity (the commutativity is understood in graded sense). Similarly, a coalgebra over ${Com}$ is a cocommutative coassociative differential graded $\kk$-coalgebra with counity.

Let us give another example of an operad. By $\Delta^k$ denote the normalized chain complex of the standard $k$-dimensional simplex. The collection $\Delta^\bullet=(\Delta^k)_{k\in\zz_+}$ forms a cosimplicial object in the category of complexes (the coface and codegeneracy operators are defined in the obvious way). For each $n\in\zz_+$, consider the tensor power $(\Delta^\bullet)^{\otimes n}=((\Delta^k)^{\otimes n})_{k\in\zz_+}$ as a cosimplicial object in the category of complexes (the coface and codegeneracy operators are defined as the tensor powers of the corresponding operators for $\Delta^\bullet$). Let $E\Delta^{(n)}=Hom(\Delta^\bullet,(\Delta^\bullet)^{\otimes n})$, where $Hom$ is taken in the category of cosimplicial objects, i.e., it is the dg-$\kk$-module whose elements of degree $d$ are all collections of maps $\Delta^k_m\to (\Delta^k)^{\otimes n}_{m+d}$ that commute with coface and codegeneracy operators. The symmetric group action on $E\Delta^{(n)}$ and the operad structure on $E\Delta$ are naturally defined. It is easy to verify (see \cite{S2}, Chapter 5, Section 2) that the chain complex of any simplicial set has a natural structure of coalgebra over $E\Delta$. In fact, this structure is universal: for any operad $A$, any natural $A$-coalgebra structure on the chain complex of a simplicial set is induced by the natural $E\Delta$-coalgebra structure via certain morphism of operads $A\to E\Delta$ (which is determined uniquely by this condition).

For any operad $A$, we define its suboperad $A^+$ by setting $(A^+)^{(0)}=0$ and $(A^+)^{(n)}=A^{(n)}$ for $n>0$ (the composition and the unit element remain the same as in $A$). We will call $A^+$ the \emph{positive part} of $A$. For example, an algebra over $Com^+$ is a commutative associative differential  $\kk$-algebra  without unit.

There exists a parallel (and simpler) theory of non-symmetric operads (see, e.g., \cite{gr_loday_vallette}, Section 5.9). In the non-symmetric case we deal with families of complexes rather than symmetric families (so, we have no action of symmetric groups). The definitions of plethysm and so on are changed accordingly. The non-symmetric operad $As$ is defined in the same way as the symmetric operad $Com$: each component $As^{(n)}$ is a rank $1$ free $\kk$-module with obvious compositions. An algebra over $As$ is an associative algebra with unit, an algebra over $As^+$ is an associative algebra without unit.

\section{$E_\infty$-operads and the category of $E_\infty$-coalgebras}

For any operad $A$ in the category of complexes, its homology $HA$ is an operad in the category of graded $\kk$-modules, and any morphism of operads induces a morphism of homology operads. We say that a morphism of operads is a \emph{quasi-isomorphism} if it induces an isomorphism in homology. Note also that any operad in the category of complexes can be viewed as an operad in the category of graded modules (by ignoring the differential).

\defn{We say that an operad $A$ in the category of graded $\kk$-modules is \emph{free} if there is a set of its elements $S=\{s_\alpha\}$: $s_\alpha\in A_{k_\alpha}^{(n_\alpha)}$, such that the following condition holds: for any operad $B$, any map $\vphi:S\to B$ such that $\vphi(s_\alpha)\in B_{k_\alpha}^{(n_\alpha)}$ has a unique extension $\tilde\vphi:A\to B$ that is a morphism of operads. An operad in the category of complexes is said to be \emph{quasifree} if it is free as an operad in the category of graded $\kk$-modules.}

The set $S$ in the definition of a quasifree operad is called the \emph{set of generators} of the operad. It is easy to see that if $A$ a quasifree operad with the set of generators $S$, then, for each $n$, there is a basis of $A^{(n)}$ indexed by pairs $(t,\sigma)$, where $t$ is an $n$-exit tree composed of vertices of types $s_\alpha$, and $\sigma\in\Sigma_n$: the corresponding element of the basis is obtained by taking compositions of the generators $s_\alpha$ in the order prescribed by the tree $t$ and then applying the permutation $\sigma$.

Suppose that all generators of a quasifree operad $A$ belong to the components $A^{(n)}$ with $n>1$. Then the differential of each generator $\partial s_\alpha=r_\alpha$ can be expressed in those generators $s_\beta$ for which $n_\beta<n_\alpha$ or $n_\beta=n_\alpha$, $k_\beta<k_\alpha$. Moreover, for any set of generators $s_\alpha$ with prescribed $k_\alpha,n_\alpha>1$, any set of such expressions for $r_\alpha$ satisfying the condition $\partial r_\alpha=0$ (to compute $\partial r_\alpha$, we use $\partial s_\beta=r_\beta$ and the Leibnitz rule) uniquely determines a quasifree operad.

Note that according to our definition, if $A$ is a quasifree operad , then $A^{(n)}$ is a complex of free $\kk\Sigma_n$-modules for each $n$.

\defn{Let $F$ be a left (right, two-sided) module over an operad $A$ in the category of graded $\kk$-modules. The module $F$ is said to be \emph{free} if it contains a set of elements $Q=\{q_\alpha\}$, $q_\alpha\in F_{k_\alpha}^{(n_\alpha)}$, satisfying the following condition: for any left (right, two-sided) $A$-module $G$, any map $\vphi:Q\to G$ such that $\vphi(q_\alpha)\in G_{k_\alpha}^{(n_\alpha)}$ extends uniquely to a morphism of left (right, two-sided) $A$-modules. A module over an operad in the category of complexes is said to be \emph{quasifree} if it becomes free after forgetting the differential.}

Suppose that $A$ is a quasifree operad and $F$ is a quasifree bimodule over it. Assume that all generators $s_\alpha$ of the operad $A$ belong to the components $A^{(n)}$ with $n>1$, and all generators $q_\beta$ of the bimodule $F$ belong to the components $F^{(n)}$ with $n>0$. Then, for each $n$, there is a basis of $F^{(n)}$ indexed by pairs $(t,\sigma)$, where $\sigma\in\Sigma_n$ and $t$ is an $n$-exit tree composed of vertices of types $s_\alpha$ and $q_\beta$ such that any monotone path along $t$ from its entry to any exit passes a vertex corresponding to a generator of $F$ exactly once. As above, the element $\sigma\in\Sigma_n$ is applied to the exit of the tree. Similarly to the case of a quasifree operad, a quasifree module can be defined recursively by expressing the differential of each generator in terms of the generators from the previous components and the generators from the same component of smaller degree.

Now we proceed to $E_\infty$-operads. We start with the standard definition.

\defn{An operad $E$ is called an \emph{$E_\infty$-operad} if all $E^{(n)}$ are free $\kk\Sigma_n$-modules and there is a quasi-isomorphism $E\to Com$.}

We will need a more narrow class of operads.

\defn{An operad $E$ is called an $E^+_\infty$-\emph{operad} if the following conditions hold:
\begin{enumerate}
\item the operad $E$ is quasifree;
\item exactly one of its generators $m\in E^{(2)}$ has degree $0$;
\item the remaining generators have positive degree and belong to $E^{(n)}$ with $n>1$;
\item there is a quasi-isomorphism $E\to Com^+$ sending $m$ to $e^{(2)}$.
\end{enumerate}}

In particular, an $E^+_\infty$-operad is a quasifree resolution  of the positive part of the commutative operad. We will give a construction of a $E^+_\infty$-operad (see \cite{S}, Chapter 5, Section 1). 
By $E(0)$ denote the free operad generated by one element $m\in E(0)^{(2)}$. Let us define a surjective homomorphism $\vphi:E(0)\to Com$ by setting $\vphi(m)=e^{(2)}\in Com(2)$. Let us chose a set of generators $v_i$ for the kernel of $\vphi$. To each of $v_i$, we assign a new free generator $w_i$ (of degree 1), set $\partial w_i=v_i$, and denote the quasifree operad generated by $m$ and all $w_i$ by $E(1)$. It is clear that the natural projection $E^{(1)}\to Com^+$ induces the isomorphism of homology in degree $0$. Further, as the $k$-th step we construct an operad $E(k)$ with a projection $Com^+$ inducing an isomorphism of homology in degree less than $k$. To this end, we add to the set of generators of $E(k-1)$ new free generators in degree $k$ so that their images under the action of the differential $\partial$ generate the $(k-1)$-th homology of the operad $E(k-1)$. The union of all operads $E(k)$ constitute the required operad $E$.

Note that any $E^+_\infty$-operad $E$ can be obtained by this construction. Indeed, for each $k$, $n$, consider the suboperad $E(k,n)\subset E$ generated by all generators of $E$ belonging to $E_{k'}^{(n')}$ with $k'\le k$ and $n'\le n$. For any generator $w\in E_k^{(n)}$, its differential $\partial w$ belongs to $E(k-1,n)$, which allows us to use induction by $k$ and $n$ in proofs and constructions. Let $E(k)=\bigcup_{n}E(k,n)$. It is clear that the sequence $E(k)$ satisfies the above description.

Let $E(\infty,n)=\bigcup_{k}E(k,n)$. (Operads of this type can be constructed directly, without reference to the enveloping $E^+_\infty$-operad.)

In the world of non-symmetric operads, we have similar notion of an $A^+_\infty$-operad, i.e., a quasifree resolution of operad $As$. In fact, such an operad $A$ can be explicitly constructed (see, e.g., \cite{Smith2}, Section 4). For each $n\ge2$, there is a unique generator $d^{(n)}\in A^{(n)}$ of degree $n-2$, and $\partial d^{(n)}=\sum_{k=2}^{n-1}\sum_{l=1}^{n-k+1}(-1)^{(k+1)(n-k+l)}d^{(k)}\circ_ld^{(n-k+1)}$.

\thrm{Let $n\ge2$. Consider the complex $C$ such that $C_k=E(\infty,n-1)^{(n)}_k$ for $k\ge0$, $C_{-1}=\kk$, and the differential $C_0\to C_{-1}$ is the projection $E(\infty,n-1)^{(n)}_0\to Com^{(n)}=\kk$. Then the only nonzero homology of $C$ is $H_{n-2}C=(Lie^{(n)})^*$, where $Lie$ is the Lie operad.}

\prf{This follows from the Koszul duality between $Com$ and $Lie$ over~$\kk$, see, e.g.,~\cite{GK,Fr}.}

\cnl{\label{min_operad}Suppose that, for each $n>1$, we are given a free $\kk\Sigma_n$-resolution of the $\kk\Sigma_n$-module $(Lie^{(n)})^*$ (i.e., a complex  $\dots L^{(n)}_2\to L^{(n)}_1\to L^{(n)}_0\to 0$ with a quasi-isomorphism to $(Lie^{(n)})^*$ in degree $0$). Let $r_k^{(n)}$ be the number of free generators of $L^{(n)}_k$. Then there exists an  $E^+_\infty$-operad with generators $m_\alpha\in E_{k_\alpha}^{(n_\alpha)}$, where the number of generators in $E_{k}^{(n)}$ is equal to $r_{k-n+2}^{(n)}$, $n\ge2$, $k\ge n-2$.}

\exmp{For $n=2$, we have the following resolution of $(Lie^{(2)})^*$ consisting of rank $1$ free $\kk\Sigma_2$-modules: \ $\dots\to\kk\Sigma_2s_2\to\kk\Sigma_2s_1\to\kk\Sigma_2s_0\to0$, where $\partial s_{0}=0$, $\partial s_{k}=s_{k-1}+(-1)^k\sigma s_{k-1}$ ($k>0$), $\sigma\in \Sigma_2$ is the transposition $(1,2)$. Thus we obtain the corresponding basis for $E^{(2)}$: we have one free generator $m^{(2)}_k$ in each degree $k\ge0$, and $(m^{(2)}_k,\sigma m^{(2)}_k)$ is a basis of $E^{(2)}_k$.

For $n=3,4$, there are also resolutions of $(Lie^{(n)})^*$ consisting of rank $1$ free $\kk\Sigma_n$-modules (we skip the details here). Thus we obtain free generators $m^{(3)}_k$, $k\ge1$, and $m^{(4)}_k$, $k\ge2$. For an arbitrary $n>4$, it is clear that the free $\kk\Sigma_n$-module $L^{(n)}_0$ in the resolution $\dots L^{(n)}_2\to L^{(n)}_1\to L^{(n)}_0\to 0$ of  $(Lie^{(n)})^*$ can be taken to be of rank $1$. Thus we can assume that in $E^{(n)}_k$ we have exactly one generator for $2\le n\le4$ and $k\ge n-2$ or $n>4$ and $k=n-2$, no generators for $k<n-2$, and finitely many generators for $n>4$, $k>n-2$. Besides, it is readily seen that the differentials for $m^{(n)}_{n-2}$ can be defined by the same formula as for the generators of the standard $A_\infty^+$-operad. Thus we can assume that our $E^+_\infty$-operad, regarded as a non-symmetric operad, contains the standard $A_\infty^+$ operad.
}

\cn{For any $n>1$, the $\kk\Sigma_n$-module $(Lie^{(n)})^*$ has a free $\kk\Sigma_n$-resolution made of rank $1$ free $\kk\Sigma_n$-modules.}

If this conjecture is true, then the following conjecture is also true.

\cn{There is a $E^+_\infty$-operad with exactly one generator in each $E^{(n)}_k$ for $n\ge2$ and $k\ge n-2$.}

The main property of $E^+_\infty$-operads is a direct consequence of the definition.

\thrm{\label{univers}Suppose that $E$ is an $E^+_\infty$-operad and $A$ is arbitrary operad such that $H_kA=0$ for $k>0$. Then, any morphism of operads $Com^+\to HA$ can be lifted to a morphism $E\to A$, and this lifting is determined uniquely up to a chain homotopy.}

In the categories of $\kk$-modules and graded $\kk$-modules, it is clear how to define the operad by generators and relations (in a similar way to defining a free operad). We define the \emph{unital binary operad} $Bin$ in the category of $\kk$-modules by generators $m\in Bin^{(2)}$ and $p\in Bin^{(0)}$ and relations $p\circ_1m=1$, $p\circ_2m=1$. Clearly, we have $Bin^{(0)}=\kk p$, and $Bin^+$ is a free operad with one generator $m$. An algebra over $Bin$ is a $\kk$-module with a bilinear binary operation (which is not assumed to be either associative or commutative) and with a neutral element for this operation. As usual, we consider $\kk$-modules as complexes concentrated in degree $0$; thus we can regard $Bin$ as an operad in the category of complexes.

\defn{An operad $E$ is called an $E^{su}_\infty$-\emph{operad} (strictly unital $E_\infty$-operad) if the following conditions hold:
\begin{enumerate}
\item all of $E^{(n)}$ are graded with nonnegative integers;
\item the operad $Bin$ is embedded in $E$ as the suboperad consisting of all elements of $E$ of degree $0$;
\item the subcomplex $\check E^{(n)}\subset E^{(n)}$ consisting of all elements $u\in E^{(n)}$ such that $p\circ_i u=0$ for $i=1,\dots,n$ is acyclic;
\item the operad $E^+$ is a $E^+_\infty$-operad with $m\in Bin^{(2)}$ serving as one of its free generators, and for all other free generators $m_\alpha\in E^{(n_\alpha)}_{k_\alpha}$ ($k_\alpha>0$), we have $m_\alpha\in\check E^{(n_\alpha)}$;
\item there is a quasi-isomorphism $E\to Com$, sending $m$ to $e^{(2)}$ and $p$ to $e^{(0)}$.
\end{enumerate}}

\thrm{Let $E$ be an $E^{su}_\infty$-operad, and let $A$ be an operad. Suppose that we are given any pair of morphisms of operads $\vphi:Com\to HA$ and $\psi_0:Bin\to A$ such that $\psi_0$ induces $\vphi$ via the morphism $Bin\cong HBin\to HE\cong Com$ induced by inclusion $Bin\to E$. Let $\check A^{(n)}=\{a\in A^{(n)}\mid \psi(p)\circ_i a=0, i=1,\dots,n\}$. If $H_k\check A^{(n)}=0$ for all $k>0$ and $n>0$, then there exists a morphism $\psi:E\to A$ coinciding with $\psi_0$ on $Bin$ and thus inducing $\vphi$. Besides, the morphism $\psi$ is determined uniquely up to a chain homotopy.
}

For example, consider the operad $A=E\Delta$. We have $H_kA=0$ for $k>0$ and $H_0A\cong Com$. Since $A^{(0)}\cong\kk$ and $A^{(1)}\cong\kk$, any choice of $e\in A^{(2)}$ that is projected to $e^{(2)}\in Com^{(2)}$ determines a morphism from $Bin$ to $A$ compatible to the isomorphism $Com\cong H_0A$ (the generator $p\in Bin^{(0)}$ is automatically projected to the basis element of $A^{(0)}\cong\kk$). It is easy to check that the complex $\check A^{(n)}$ is acyclic for any $n>0$. Thus, for any $E^{su}_\infty$-operad $E$ we get a morphism of $E$ to $A=E\Delta$ determined uniquely up to a chain homotopy. Thus the chain complex of any simplicial set has a natural structure of $E$-coalgebra, and this natural structure is determined uniquely up to a natural chain homotopy.

Let us choose an $E^{su}_\infty$-operad $E$. We will consider $E$-coalgebras which are complexes of free $\kk$-modules. By a homomorphism of $E$-coalgebras $f:X\to Y$ we mean a morphism of complexes which commutes with the action of $E$. This condition can be described as follows. We have the homomorphisms of $E$ to $\cee^X$ and to $\cee^Y$ that define the $E$-coalgebra structures on $X$ and $Y$. Since $\cff^{X,Y}$ is a left $\cee^Y$-module and a right $\cee^X$-module, we obtain a structure of $E$-bimodule on $\cff^{X,Y}$. The operad $E$ can also be regarded as a bimodule over itself. Note that a morphism of complexes $f:X\to Y$ is a homomorphism of $E$-coalgebras if and only if there exists a homomorphism of  $E$-bimodules $E\to\cff^{X,Y}$ sending $1\in E^{(1)}$ to $f\in(\cff^{X,Y})^{(1)}$ (this condition determines such a homomorphism uniquely). Thus the set of homomorphisms of $E$-coalgebras $X\to Y$ can be identified with the set of homomorphisms of $E$-bimodule $E\to\cff^{X,Y}$.

Consider the operad $Com$ as a bimodule over the operad $E$ (via the projection $E\to Com$).

\prop{\label{propf}There exists a non-negatively graded $E$-bimodule $F$ such that
\begin{enumerate}
\item all of $F^{(n)}$ are graded with nonnegative integers;
\item the components $F^{(0)}$ and $F^{(1)}$ are free $\kk$-modules of rank $1$ with generators $f^{(0)}$ and $f^{(1)}$ (of degree $0$), respectively, which satisfy the relations $(p)\circ f^{(1)}=f^{(0)}=()\circ p$;
\item the subcomplex $\check F^{(n)}\subset F^{(n)}$ consisting of all elements $w\in F^{(n)}$ such that $p\circ_i w=0$ for $i=1,\dots,n$ is acyclic;
\item $F^+$ is a quasifree bimodule over the operad $E^+$, with $f^{(1)}$ serving as one of its free generators, and for all other free generators $f_\beta\in F^{(n_\beta)}_{k_\beta}$, we have $k_\beta>0$ and $f_\beta\in\check F^{(n_\beta)}$;
\item there is a quasi-isomorphism of $E$-bimodules $F\to Com$, sending $f^{(1)}$ to $1=e^{(1)}$ and $f^{(0)}$ to $e^{(0)}$.
\end{enumerate}
Moreover, suppose that the quasifree operad $E^+$ is generated by a collection of free generators $(m_\alpha)$, where $m_\alpha\in E^{(n_\alpha)}_{k_\alpha}$. Than the bimodule $F$ can be chosen in such a way that the $E^+$-bimodule $F^+$ a collection of free generators consisting of one generator $f_0$ in $F^{(1)}_0$ and certain generators $f_\alpha\in F^{(n_\alpha)}_{k_\alpha+1}$ that are in one-to-one correspondence with the generators of $E^+$ (this correspondence preserves arity $n$ and increases degree $k$ by $1$).
}

\prf{
As in the case of $E^+_\infty$-operads, we use induction by $k$. The differential of each generator $\nu_\alpha\in F^{(n_\alpha)}_{k_\alpha+1}$ is expressed in the form $\partial f_\alpha=m_\alpha\circ f_0-(f_0,\dots,f_0)\circ m_\alpha+\dots$, where dots denote some expression in $m_\beta$ and $f_\beta$ with $k_\beta<k_\alpha$.
}

Let us call any $E$-bimodule $F$ satisfying the conditions stated in Proposition~\ref{propf} a \emph{universal bimodule} over $E$. This definition is motivated by the following properties of $F$.

\prop{\label{propuni1}
Let $B$ be an $E^+$-bimodule such that $H_kB=0$ for $k>0$. Suppose we are given a homomorphism of $Com^+$-bimodules $\vphi:Com^+\to H_0B$. Then there exists a homomorphism of $E^+$-bimodules $\psi:F^+\to B$ that induces the homomorphism $\vphi$ in homology, and $\psi$ is determined uniquely up to a chain homotopy.
}

\prop{\label{propuni2}
Suppose that $B$ is an $E$-bimodule and $g^{(0)}\in B^{(0)}_0$, $g^{(1)}\in B^{(1)}_0$ satisfy the relations $(p)\circ g^{(1)}=g^{(0)}=()\circ p$. Let $\check B^{(n)}=\{b\in B^{(n)}\mid p\circ_i b=0, i=1,\dots,n\}$. Suppose that $H_k\check B^{(n)}=0$ for all $k>0$ and $n>0$. If we are given a homomorphism of $Com$-bimodules $\vphi:Com\to H_0B$ such that $\vphi(e^{(0)})=\bar g^{(0)}$, $\vphi(e^{(1)})=\bar g^{(1)}$, where $\bar g^{(0)},\bar g^{(1)}$ are the images of $g^{(0)},g^{(1)}$ in $H_0B$, then there exists a homomorphism of $E$-bimodules $\psi:F\to B$ inducing the homomorphism $\vphi$ in homology and sending $f^{(0)}$ to $g^{(0)}$ and $f^{(1)}$ to $g^{(1)}$. Such a homomorphism is determined uniquely up to a chain homotopy.
}

Let us choose an $E^{su}_\infty$-operad $E$ and a universal $E$-bimodule $F$. We define the category $E$-$coalg$ as follows: its objects are  $E$-coalgebras, and the morphisms from $K$ to $L$ are classes of $E$-bimodule homomorphisms $F\to\cff_{K,L}$ up to chain homotopy. Note that if a morphism $\bar\vphi$ is represented by an $E$-bimodule homomorphism $\vphi:F\to\cff_{K,L}$, then $\vphi(f^{(1)})$ is a well-defined chain map (morphism of complexes) $K\to L$. We say that the morphism $\bar\vphi$ of the category $E$-$coalg$ \emph{extends} the chain map $\vphi(f^{(1)})$. (The definition of morphisms with the help of an operad bimodule was proposed in \cite{Lyu} for the case of $A_\infty$-algebras.)
 
Similarly, in the category of $E^+$-coalgebras the morphisms from $K$ to $L$ are classes of $E^+$-bimodule homomorphisms $F^+\to\cff_{K,L}$ up to chain homotopy. The universality of the bimodule $F$ implies that the composition of morphisms is well-defined and is associative. Similar reasoning based on universality shows that another choice of an $E^{su}_\infty$-operad and a universal bimodule gives us equivalent category of coalgebras.

\section{The cobar construction in the category $E$-$coalg$ and the dimension completion of the fundamental group}

The ring $\kk$, regarded as a dg-$\kk$-module (concentrated in degree $0$), has a trivial structure of $E$-coalgebra (because $\kk$ has a standard structure of $Com$-coalgebra, and we have a projection $E\to Com$). Note that, for any $E$-coalgebra $K$, the counit, i.e., the action of $p\in Bin^{(0)}=E^{(0)}$, defines a homomorphism of $E$-coalgebras $C\to\kk$. Denote this homomorphism by $\pi$. An $E$-coalgebra $K$ is said to be  \emph{coaugmented} if it is equipped with a homomorphism of $E$-coalgebras $\iota:\kk\to K$ right inverse to $\pi$ (that is, $\pi\circ\iota=\id_\kk$). If $K$ is a chain complex of a simplicial set, and the structure of $E$-coalgebra is defined as above, then choosing a coaugmentation is equivalent to choosing a base point---one of the vertices of our simplicial set.

\prop{\label{propaug}Let $K$ be a coaugmented $E$-coalgebra. For each of the free generators $m_\alpha \in E^{(n_\alpha)}_{k_\alpha}$ of $E^+$, we set $\tilde m_\alpha=((\id_K-\iota\circ\pi)\otimes\dots\otimes(\id_K-\iota\circ\pi))\circ m_\alpha$. Then the map $m_\alpha\mapsto\tilde m_\alpha$ endows  $\tilde K=\Ker\pi\subseteq K$ with an $E^+$-coalgebra structure.
}

Consider the category whose objects are coaugmented $E$-coalgebras and whose morphisms are the morphisms in $E$-$coalg$ for which the image of $f^{(1)}$ in $\cff_{K,L}$ is a map commuting with coaugmentations. It is easy to check that the correspondence $K\mapsto \tilde K=\Ker\pi$ gives us a functor from the category of coaugmented $E$-coalgebras to the category of $E^+$-coalgebras.

Let $K$ be an $E^+$-coalgebra. Consider $\overline{TK[1]}$, i.e., the completed tensor algebra of the complex $K[1]$ (the algebra of non-commutative formal power series). Be $K[1]$ we denote the desuspension of the complex $K$ (that is, $K[1]_k=K_{k+1}$, and the differential differs from the differential in $K$ by the factor $(-1)$). Since the non-symmetric operad $A_\infty$ is embedded into $E$, we have the differential $D=\partial+\sum d^{(n)}$ on $\overline{TK[1]}$ (the cobar construction for $A_\infty$-coalgebras, see \cite{Smith2}). Consider $\overline{TK[1]}$ as a complex with respect to the differential $D$. It is endowed with augmentation $\pi:\overline{TK[1]}\to\kk$, which sends each formal power series to its constant term, and coaugmentation $\iota:\kk\to\overline{TK[1]}$, which sends each constant to a constant power series.

\thrm{The augmentation $\pi$ and the coaugmentation $\iota$ on the cobar construction $\overline{TK[1]}$ can be extended to a natural structure of coaugmented $E$-coalgebra, and this $E$-coalgebra structure is determined uniquely up to a chain homotopy. More precisely, such extension is a functor from the category of $E^+$-coalgebras to the category of coaugmented $E$-coalgebras, and this functor is unique up to isomorphism.}

Taken together with Proposition \ref{propaug}, this theorem allows to iterate the cobar construction. Another variant of iterating the cobar construction is studied in \cite{Smith2}. By the uniqueness, these constructions are equivalent. Applying the results of~\cite{Smith2} to our construction, we can prove that if we start from the chain complex of a simply connected simplicial set, then the cobar construction is isomorphic (in our category $E$-$coalg$) to the chain complex of the loop space.

However, we are mostly interested in the case of spaces that are not simply connected.

\thrm{
Suppose that the initial $E$-coalgebra $K$ is the chain complex of a simplicial set $X$ with a single vertex $x_0$. Then the $E$-coalgebra structure on the cobar construction $\overline{T\tilde K[1]}$ together with the tensor multiplication define a Hopf algebra structure on $H_0\overline{T\tilde K[1]}$. This Hopf algebra is isomorphic to the completion of the group algebra of the fundamental group $\pi_1(X,x_0)$ with respect to the powers of the augmentation ideal (see \cite{Ry2}). Thus, the closure of the image of $\pi_1(X,x_0)$ in $H_0\overline{T\tilde K[1]}$ can be obtained as the set of group-like elements of this Hopf algebra. Thus, for $\kk=\zz$, the $E$-coalgebra structure on the chain complex of a simplicial set determines the dimension completion of the fundamental group (see \cite{MP}).
}

Suppose that complexes $K$ and $L$, morphisms $f:K\to L$ and $g:L\to K$, and a map $h:K\to K$ of degree $1$ satisfy the following conditions: $fg=\id_L$, $gf=\id_K+[h,\partial]$, $h^2=0$, $fh=0$, $hg=0$. In this case we say that we are given a \emph{strong deformation retraction} of $K$ onto $L$.

\thrm{Suppose that $K$ is a $E$-coalgebra, and we are given a strong deformation retraction $(f,g,h)$ of $K$ onto $L$. Then there exists an $E$-coalgebra structure on $L$ such that $f$ and $g$ can be extended to isomorphisms in the category $E$-$coalg$.}

If $\kk=\zz$ and the homology of $K$ is a free $\zz$-module, then there is a strong deformation retraction $(f,g,h)$ of $K$ onto $HK$. Hence the homology has a natural $E$-coalgebra structure. In the case when $K$ is the chain complex of a connected simplicial set we obtain the following result.

\prop{\label{free_homology} Suppose that the integer homology groups $H_*(X)$ of a connected simplicial set $X$ are free $\zz$-modules. Then the homology has a natural $E$-coalgebra structure, which is unique up to an isomorphism in the category $E$-$coalg$, and this structure determines the dimension completion of the fundamental group of $X$.}

\exmp{Let us find some invariants of $E$-coalgebra structures on $H_*(X)$ with respect to isomorphisms in the category $E$-$coalg$ which extend the identity map. We assume that the operad $E$ is constructed as in Corollary~\ref{min_operad} and the universal bimodule $F$ satisfies the properties stated in Proposition~\ref{propf}. Thus, for small $n$, in $E$ we have the generators $m_0^{(0)}=p$, $m_0^{(1)}$ (the unit element),  $m_k^{(2)}$ for each $k\ge0$, and $m_k^{(3)}$ for each $k\ge1$. In $F$ we have the generators $f_0^{(0)}=()\circ p$, $f_0^{(1)}$, $f_k^{(2)}$ for each $k\ge1$, and $f_k^{(3)}$ for each $k\ge2$.

For each of these generators $m_k^{(n)}$, by $\hat{m}_k^{(n)}$ and $\hat{\hat{m}}_k^{(n)}$ we denote the corresponding operators $H_*(X)\to H_*(X)^{{\otimes}n}$ according to two $E$-coalgebra structures on $H_*(X)$. The operators $\hat{f}_k^{(n)}:H_*(X)\to H_*(X)^{{\otimes}n}$ represent a $E$-coalgebra morphism that extends the identity map between these $E$-coalgebras. Thus $\hat{f}_0^{(1)}=\id=\hat{m}_0^{(1)}=\hat{\hat{m}}_0^{(1)}$.

Let $\sigma$ be the transposition $(1,2)\in\Sigma_2$. We have the following relations:
\begin{align*}
\partial m^{(2)}_{k+1}&=m^{(2)}_k-(-1)^k\sigma m^{(2)}_k,\quad k\ge0;\\
\partial f^{(2)}_{1}&=m^{(2)}_0\circ f^{(1)}_0-(f^{(1)}_0,f^{(1)}_0)\circ m^{(2)}_0;\\
\partial f^{(2)}_{k+1}&=m^{(2)}_k\circ f^{(1)}_0-(f^{(1)}_0,f^{(1)}_0)\circ m^{(2)}_k-(f^{(2)}_{k}+(-1)^k\sigma f^{(2)}_{k}),\quad k\ge1;\\
\partial m^{(3)}_{1}&=m^{(2)}_0\circ_1m^{(2)}_0 - m^{(2)}_0\circ_2m^{(2)}_0;\\
\partial f^{(3)}_{2}&=m^{(3)}_1\circ f^{(1)}_0-(f^{(1)}_0,f^{(1)}_0,f^{(1)}_0)\circ m^{(3)}_1 -\\ &-m^{(2)}_0\circ_1f^{(2)}_{1}+m^{(2)}_0\circ_2f^{(2)}_{1}-f^{(2)}_{1}\circ_1m^{(2)}_0+f^{(2)}_{1}\circ_2m^{(2)}_0.
\end{align*}
Since the differential on homology is equal to zero, it follows that
\begin{align*}
\hat{m}^{(2)}_k&=(-1)^k\sigma \hat{m}^{(2)}_k,\quad k\ge0;\\
\hat{\hat{m}}^{(2)}_k&=(-1)^k\sigma \hat{\hat{m}}^{(2)}_k,\quad k\ge0;\\
\hat{\hat{m}}^{(2)}_0&=\hat{m}^{(2)}_0;\\
\hat{\hat{m}}^{(2)}_k-\hat{m}^{(2)}_k&=\hat{f}^{(2)}_{k}+(-1)^k\sigma \hat{f}^{(2)}_{k},\quad k\ge1;\\
\hat{m}^{(2)}_0\circ_1\hat{m}^{(2)}_0 &= \hat{m}^{(2)}_0\circ_2\hat{m}^{(2)}_0;\\
\hat{\hat{m}}^{(2)}_0\circ_1\hat{\hat{m}}^{(2)}_0 &= \hat{\hat{m}}^{(2)}_0\circ_2\hat{\hat{m}}^{(2)}_0;\\
\hat{\hat{m}}^{(3)}_1-\hat{m}^{(3)}_1&=\hat{\hat{m}}^{(2)}_0\circ_1\hat{f}^{(2)}_{1}+\hat{\hat{m}}^{(2)}_0\circ_2\hat{f}^{(2)}_{1}
-\hat{f}^{(2)}_{1}\circ_1\hat{m}^{(2)}_0+\hat{f}^{(2)}_{1}\circ_2\hat{m}^{(2)}_0.
\end{align*}

From these equation we see, of course, that the comultiplication in homology is cocommutative and coassociative. We also see that it is well-defined, i.e., invariant with respect to isomorphisms in the category $E$-$coalg$.

In view of applications to the fundamental group of $X$, let us consider those components of the above operators that touch only $H_1=H_1(X)$ and $H_2=H_2(X)$. Taking into account the Koszul sign convention and the fact that the degree of $H_1$ is odd and the degree of $H_2$ is even, we see that after forgetting this grading we can regard $\hat{m}^{(2)}_0=\hat{\hat{m}}^{(2)}_0$ as a map from $H_2$ to $\Lambda^2H_1\cong[H_1,H_1]$ (the commutator in the free associative algebra generated by $H_1$). Next, the operator $\hat{m}^{(2)}_1$ has components $H_1\to H_1\otimes H_1$ and $H_2\to (H_1\otimes H_2\oplus H_2\otimes H_1)$ (with image in symmetric tensors). We see that it is not invariant: these components can be changed by adding an arbitrary map taking values in anticommutators, and can be changed by adding an arbitrary map taking values in commutators. Thus the component $H_2\to (H_1\otimes H_2\oplus H_2\otimes H_1)$ can be made arbitrary, while the component $H_1\to H_1\otimes H_1$ has an invariant which is a dualized form of the Steenrod square operation in cohomology.

Suppose that the $\hat{m}^{(2)}_0=\hat{\hat{m}}^{(2)}_0$, $\hat{m}^{(2)}_1$, and $\hat{\hat{m}}^{(2)}_1$ are given so that $\hat{f}^{(2)}_{1}$ can be chosen to satisfy the above relation. Then we have additional freedom: we can add any operator $H_1\to [H_1,H_1]$ to $\hat{f}^{(2)}_{1}:H_1\to H_1\otimes H_1$ and any operator $H_2\to [H_1,H_2]$ to  $\hat{f}^{(2)}_{1}:H_2\to(H_1\otimes H_2\oplus H_2\otimes H_1)$. Thus we get an invariant element in the Abelian group
$$
\frac{Hom(H_2,[H_1,[H_1,H_1]]/[H_1,\hat{m}^{(2)}_0(H_1)])}{\delta Hom(H_1,[H_1,H_1])}\;,
$$
where the map $\delta:Hom(H_1,[H_1,H_1])\to Hom(H_2,[H_1,[H_1,H_1]]/[H_1,\hat{m}^{(2)}_0(H_1)])$ sends each $\mu\in Hom(H_1,[H_1,H_1])$ to the class of $\mu\circ_1\hat{m}^{(2)}_0-\mu\circ_2\hat{m}^{(2)}_0$. It is readily seen that this invariant is a dualized form of the triple Massey product on $H^1$ (see~\cite{Ry2}).

Let $X$ be the simplicial set of singular simplexes of the complement to a complex hyperplane arrangement. Then the assumptions of Proposition~\ref{free_homology} are satisfied. Moreover, the map $\hat{m}^{(2)}_0:H_2\to [H_1,H_1]$ is injective. It is shown in \cite{Ry2} that the the invariant described above is an invariant of the fundamental group with respect to isomorphisms that induce the canonical isomorphism of the first homology groups. This invariant was first constructed in \cite{Ry} (by means of the lower central series of the fundamental group). It was used there to prove that there exist an example of combinatorially equivalent arrangements with non-isomorphic fundamental groups.
}

\end{document}